\documentclass[conference]{IEEEtran}
\IEEEoverridecommandlockouts
\usepackage{cite}
\usepackage{amsmath,amssymb,amsfonts}
\usepackage{algorithmic}
\usepackage{graphicx}
\usepackage{textcomp}
\usepackage{multirow}
\usepackage{float}
\usepackage{subcaption}
\usepackage{xcolor}
\def\BibTeX{{\rm B\kern-.05em{\sc i\kern-.025em b}\kern-.08em
    T\kern-.1667em\lower.7ex\hbox{E}\kern-.125emX}}
\let\vec\mathbf
\newcommand{\T}{^\mathsf{T}}
\begin{document}

\title{Homotopy Continuation for Sensor Networks Self-Calibration}

\author{\IEEEauthorblockN{Luca Ferranti}
\IEEEauthorblockA{
\textit{University of Vaasa}\\
Vaasa, Finland}
\and
\IEEEauthorblockN{Kalle Åström}
\IEEEauthorblockA{
\textit{Lund University}\\
Lund, Sweden}
\and
\IEEEauthorblockN{Magnus Oskarsson}
\IEEEauthorblockA{
\textit{Lund University}\\
Lund, Sweden}
\and
\IEEEauthorblockN{Jani Boutellier}
\IEEEauthorblockA{\textit{University of Vaasa}\\
Vaasa, Finland}
\and
\IEEEauthorblockN{Juho Kannala}
\IEEEauthorblockA{
\textit{Aalto University}\\
Espoo, Finland}
}

\maketitle

\begin{abstract}
Given a sensor network, TDOA self-calibration aims at simultaneously estimating the positions of receivers and transmitters, and transmitters time offsets. This can be formulated as a system of polynomial equations. Due to the elevated number of unknowns and the nonlinearity of the problem, obtaining an accurate solution efficiently is nontrivial. Previous work has shown that iterative algorithms are sensitive to initialization and little noise can lead to failure in convergence. Hence, research has focused on algebraic techniques. Stable and efficient algebraic solvers have been proposed for some network configurations, but they do not work for smaller networks. In this paper, we use homotopy continuation to solve four previously unsolved configurations in 2D TDOA self-calibration, including a minimal one. As a theoretical contribution, we investigate the number of solutions of the new minimal configuration, showing this is much lower than previous estimates. As a more practical contribution, we also present new subminimal solvers, which can be used to achieve unique accurate solutions in previously unsolvable configurations. We demonstrate our solvers are stable both with clean and noisy data, even without nonlinear refinement afterwards. Moreover, we demonstrate the suitability of homotopy continuation for sensor network calibration problems, opening prospects to new applications.
\end{abstract}


\section{Introduction}
\label{ch:intro}
Given a network of receivers and transmitters, network self-calibration refers to the simultaneous estimation of positions of both receivers and transmitters \cite{miluzzo2008, wendberg2013}. Network self-calibration is essential in several applications, such as beamforming \cite{ochiai2005} or structure from sound \cite{thrun2005}. The case of synchronized networks, where all time delays are known and only the positions of the nodes need to be estimated, has been greatly investigated in the literature \cite{crocco2011, kuang2013, raykar2005, kuang2013} and it is mainly solved. The case of \textit{unsynchronized} networks, where also some time offsets of the nodes need to be estimated, has proven more challenging and is still an active research area \cite{burgess2012, moses2003}. In this paper, we focus on 2D self-calibration of networks with synchronized receivers but unsynchronized transmitters, known as network 2D Time Difference Of Arrival (TDOA) self-calibration \cite{kuang2013tdoa, pollefeys2008}. 

Let a network with $m$ receivers and $n$ transmitters, from now on shortly denoted as $mr/ns$, the TDOA self-calibration can be formalized by a set of $mn$ equations in the form
\begin{equation}
    \Vert\vec{r}_i-\vec{s}_j\Vert^2=(f_{ij}-o_j)^2=d_{ij}^2,~i=1,\ldots,m~j=1,\ldots,n \label{eq:primal}
\end{equation}
where $\vec{r}_i,\vec{s}_j$ are the unknown positions of receiver $i$ and transmitter $j$, $d_{ij}$ is the distance between them, $o_j$ is the unknown transmitter offset and $f_{ij}$ is the \textit{pseudorange}, i.e. what the receiver $i$ measures from the transmitter $j$. To study the solvability of this problem, the \textit{excess constraint} is defined as the number of constraints minus the number of unknowns. For 2D TDOA, the excess constraint is
\begin{equation}
    mn-2m-3n+3, \label{eq:dof}
\end{equation}
where the term $+3$, called \textit{Gauge freedom}, is introduced because only distances are measured and hence the positions can be recovered only up to Euclidean transformation. If the expression in \eqref{eq:dof} is negative, then we have less constraints than unknowns, and the problem will have infinitely many solutions, if it is zero, then the problem is determined and it will have a finite, but not necessarily unique, number of finite solutions. These determined problems are referred to as \textit{minimal}, because they are the smallest problems to have a finite solution. Finally, if the expression in \eqref{eq:dof} is positive, then the problem is overdetermined, or \textit{subminimal}, and it will have a unique solution (given zero noise there will a unique solution that satisfies the constraints exactly, otherwise we find a solution that minimizes the constraint equations).

Given the strong nonlinearity and the high number of unknowns, numerical iterative approaches are very sensitive to initialization and hence perform poorly with a random or grid-searched initial value \cite{biswas2004, priyantha2003}. Algebraic techniques have proved more efficient for both synchronized \cite{kuang2013, larsson2020} and unsynchronized \cite{kuang2013tdoa, batstone2019, ferranti2020, pollefeys2008} network calibration. These algebraic techniques, based on Gröbner bases, are split into two phases \cite{kukelova2008}: an offline phase, during which an optimized solver for a specific kind of problem is generated; and an online phase, during which the generated solver can be applied to quickly solve any instance of that problem. While these techniques have proved successful in several applications \cite{larsson2017, kukelova2011}, they are somehow bounded by the complexity of the problem and they have not been able to solve more challenging TDOA configurations. \textit{Homotopy continuation}, a numerical algorithm for systems of polynomial equations \cite{morgan1987}, has on the other hand been shown to solve more complex problems, where Gröbner bases methods failed \cite{fabbri2020, duff2020}. In this paper, we demonstrate that homotopy continuation can also be applied to self-calibration of sensor networks, presenting stable solvers for several previously unsolvable configurations. To the best of our knowledge, we are the first to apply homotopy continuation to sensor networks self-calibration.

\begin{table}[t]
    \centering
    \caption{Different TDOA configurations for $m$ receivers and $n$ transmitters. X: previously solved. O: solved in this paper. *: reducible to a solved configuration. -: underdetermined. u: unsolved. Underlined the minimal configurations}
    \label{tab:probs}
    \begin{tabular}{cc|cccc}
        &&\multicolumn{4}{c}{$n$}\\
         &&3&4&5&6  \\\hline
         \multirow{6}{*}{$m$}&4&-&-&\underline{u}&u\\
         &5&-&O&O&X\\
         &6&\underline{O}&O/X&*&*\\
         &7&O&X&*&*\\
         &8&*&*&*&*\\
         &9&X&*&*&*\\
    \end{tabular}
\end{table}

The contribution of our paper is mainly two-fold. First, we show that homotopy continuation can robustly solve sensor network localization problems, opening frontiers for even more applications in this domain. Secondly, through homotopy continuation we present several new algebraic solvers for network configurations that were previously unsolvable. Particularly, we present solvers for the following previously unsolvable configurations (also summarized in Table~\ref{tab:probs}).
\begin{itemize}
    \item \textbf{6r/3s}: this is a minimal configuration, i.e. it is the smallest network to have a finite number of solutions. The solution, however, is not unique. From a theoretical perspective, it is interesting to investigate how many solutions the problem can have. It was speculated in \cite{stewenius2005}, that this configuration could have up to 150 distinct solutions, but some of those might be false roots. Using our homotopy solver, we are able to give empirical evidence that the number of solutions is actually much smaller.
    \item \textbf{6r/4s} and \textbf{7r/3s}: these can be regarded as the subminimal configurations of the previous ones, obtained adding one extra point to achieve uniqueness of the solution. While the 7r/3s solver is new, an algebraic solver for 6r/4s was already presented in \cite{ferranti2020}. This, however, required an extra nonlinear refinement of the solution, while our homotopy continuation solver is already stable without this post-processing step.
    \item \textbf{5r/4s}: This configuration has an interesting property, while it is not minimal, it is the only subminimal configuration which cannot be reduced to a minimal one.  Hence, it requires its own specialized solver.
    \item \textbf{5r/5s}: This solver can be considered the subminimal configuration of 4r/5s (which remains unsolved), with an extra point to ensure the uniqueness of the solution. Even if a solver for the 4r/5s were found, it would still have multiple solutions. This solver, hence, is the smallest solver to achieve a unique solution in a network with five transmitters.
\end{itemize}

This paper is structured as follows: The theoretical background of homotopy continuation is discussed in Section \ref{ch:homotopy} and the algorithms for network calibration are presented in Section \ref{ch:alg}. Finally, the results are discussed in Section \ref{ch:results} and conclusions are drawn in Section \ref{ch:conclusions}.

\section{Background: homotopy continuation}
\label{ch:homotopy}
Homotopy continuation is an iterative algorithm to solve systems of polynomial equations \cite{morgan1987} and it has several success stories in e.g. computer vision \cite{kileel2017, fabbri2020}, signal processing \cite{malioutov2005, fanjul2017}, and process design \cite{wayburn1987, chang1988} 

Let $F(\vec{x})=[f_1,\ldots,f_n]\T$ and $G(\vec{x})=[g_1,\ldots,g_n]\T$  be vectors of polynomials from the ring $\mathbb{C}[x_1,\ldots,x_n]$. Suppose we want to solve the system $F(\vec{x})=\vec{0}$ and that the system $G(\vec{x})=\vec{0}$, referred to as \textit{starting system}, can be solved easily and has at least as many distinct roots as $F$. We can now define the \textit{homotopy}
\begin{equation}
    H(\vec{x},t)=(1-t)F(\vec{x})+\gamma tG(\vec{x}),
\end{equation}
where $t$ is a new variable and $\gamma$ is a complex number with $\Vert\gamma\Vert=1$, introduced for numerical stability reasons \cite{morgan1987} . It is now easy to see that $H(\vec{x}, 1)=\gamma G(\vec{x})$ and $H(\vec{x}, 0)=F(\vec{x})$. Also if $\vec{x}_0$ is a root of $G$, then $H(\vec{x}_0,1)=\vec{0}$ and if $\vec{x}^*$ is a root of $F$, then $H(\vec{x}^*, 0)=\vec{0}$. 

The core idea is that if we have a solution $\vec{x}_t$ of $H(\vec{x}, t)=\vec{0}$, then for a small perturbation $h$ the solution $\vec{x}_{t+h}$ of $H(\vec{x}_{t+h},t+h)=\vec{0}$ will be close to $\vec{x}_t$, and hence it can be computed fast with Newton method using $\vec{x}_t$ as initial values. Hence, as $t$ varies smoothly from $1$ to $0$, $\vec{x}_0$ will smoothly converge to a root $\vec{x}^*$ of $F$. Obviously, this process has to be iterated for each root $\vec{x}_0$ of $G$.

It remains to decide how to choose the starting system. Suppose $F$ has $n$ polynomials of degree $d_1,\ldots,d_n$. By Bézout theorem \cite{cox98}, $F$ can have at most $K=\prod_{i=1}^n d_i$ distinct roots. Hence choosing $G$ to have exactly $K$ distinct solutions will guarantee finding all solutions of $F$. A system satisfying this property can be easily constructed by choosing
\begin{equation}
   G=\begin{bmatrix}x_1^{d_1}-a_1\\\vdots\\x_n^{d_n}-a_n\end{bmatrix},\label{eq:starting}
\end{equation}
where $a_i$ are some non-zero complex numbers. The system $G(\vec{x})=\vec{0}$ will then have exactly $K$ solutions.

Now, for each root $\vec{x}_0$ of $G$, we apply homotopy continuation and track the path to the roots of $F$. It is good to notice that $F$ needs not have $K$ distinct solutions. 
Indeed, sparse systems are very likely to have strictly less distinct solutions than $K$. Practically, this means that a starting system built with \eqref{eq:starting} will introduce some computational overhead, as it will track more paths than necessary. Recently, more efficient initializations, exploiting polyhedral geometry, have been proposed for sparse systems \cite{huber1995, jensen2016}. In our numerical experiments, we used the Julia implementation \textsf{HomotopyContinuation.jl} \cite{breiding2018n} and the polyhedral starting system proposed in \cite{jensen2016} as initialization.

\section{Solver design}
\label{ch:alg}
Directly applying homotopy continuation to the equations as in \eqref{eq:primal} would result in a poorly conditioned solver. In this section we describe the algebraic manipulation that leads to a more stable formulation of the problem.

\subsection{Parametrization of positions}
In the general 2D case, node positions require $2m+2n-3$ unknowns. However, it was shown in \cite{kuang2013}, that for the case with $n=3$, the positions can be parametrized using only five unknowns. 

We start by defining the \textit{compaction matrix} $\vec{\tilde{D}}\in\mathbb{R}^{(m-1)\times(n-1)}$ as
\begin{equation}
    [\vec{\tilde{D}}]_{ij} = d_{i+1,j+1}^2-d_{1,j+1}^2-d_{i+1,1}^2+d_{11}^2, \label{eq:compaction}
\end{equation}
and fix the Gauge freedom by setting $\vec{r}_1=\vec{0}$ and $\vec{r}_{2y}=0$. The remaining points can be parametrized as follows \cite{ferranti2020}
\begin{equation}
        \begin{split}
        &\vec{s}_1=\vec{L}\vec{b},\\
        &\vec{r}_i=\vec{L}^{-\mathsf{T}}\tilde{\vec{D}}_{i-1}\T,\quad i=2\ldots m,\\
        &[\vec{s}_2~\vec{s}_3]=\vec{L}\left(-\frac{1}{2}\vec{I}+[\vec{b}~\vec{b}]\right),
    \end{split} \label{eq:toa_upgrade}
\end{equation}
where $\vec{b}$ is a vector of $2$ unknowns and $\vec{L}$ is a $2\times2$ upper triangular matrix, $\vec{I}$ is the identity matrix and $\vec{\tilde{D}}_i$ denotes the $i$th row of $\vec{\tilde{D}}$. Defining the symmetric matrix $\vec{H}=(\vec{L}\T\vec{L})^{-1}$ we can derive the equations
\begin{equation}
    \begin{split}
        &\mathrm{(A)}\quad d_{11}^2=\vec{b}\T\vec{H}^{-1}\vec{b},\\
    &\mathrm{(B)}\quad d_{1j}^2-d_{11}^2=\frac{1}{4}\vec{e}_{j-1}\T\vec{H}^{-1}\vec{e}_{j-1}-\vec{b}\T\vec{H}^{-1}\vec{e}_{j-1},\\
    &\mathrm{(C)}\quad d_{i1}^2-d_{11}^2=\tilde{\vec{D}}_{i-1}\vec{H}\tilde{\vec{D}}_{i-1}\T-2\vec{b}\T\tilde{\vec{D}}_{i-1}\T,
    \end{split}\label{eq:dual}
\end{equation}
where $[\vec{e}_1,~\vec{e}_2]$ are the standard basis of $\mathbb{R}^2$ and $i=2,\ldots,m$ and $j=2,3$. Since $\vec{H}$ is symmetric, it has only $3$ unknowns. Recalling we also have $3$ unknown offsets, the total number of unknowns will be $8$. For a general configuration with $m$ receivers and $3$ transmitters we will have $1$ equation of type (A), $2$ equations of type (B) and $m-1$ equations of type (C), hence $m+2$ constraints in total. 

\subsection{Offset constraints}
It can be shown \cite{kuang2013tdoa}, that the compaction matrix \eqref{eq:compaction} must have rank $2$. For configurations with $m\geq4$ and $n\geq4$, this means that some extra constraints can be imposed on the offsets by setting all $3\times 3$ determinants in $\tilde{D}$ to $0$. While we can obtain $\binom{m-1}{3}\cdot\binom{n-1}{3}$ rank constraints in total, only $(m-3)(n-3)$ of those are independent. This approach was used in \cite{kuang2013} to solve some overdetermined problems by first solving for the offsets separately and then localizing the nodes using synchronized network calibration techniques. 
Note that this approach alone is not suitable for smaller networks, but it can still offer extra constraints to add to \eqref{eq:dual}.

During the experiments we will use both the initial formulation and the one described in this section. To distinguish between the two, we call the one using equations as in \eqref{eq:primal} the \textit{primal formulation}, and the one using equations \eqref{eq:dual} and eventual offset constraints the \textit{dual formulation}.

\subsection{Summing up: building the homotopy solvers}
The dual formulation significantly reduced the complexity of the problem and it can thus be efficiently solved by homotopy continuation. Here we discuss some solver-specific implementation details, to allow reproducibility of the experiments.
\begin{itemize}
    \item \textbf{6r/3s}: This is a minimal configuration. Using the dual formulation, we obtain a total of $8$ unknowns and $6+2=8$ equations like \eqref{eq:dual}. This nonlinear polynomial system can now be solved with homotopy continuation. Next $\vec{L}$ can be retrieved with Cholesky factorization and the positions can finally be computed with \eqref{eq:toa_upgrade}. Being a nonlinear minimal problem, the uniqueness is not guaranteed. However, not all the solutions of the dual formulation are necessarily solutions of the primal. Particularly, we accept only \textit{real solutions} of the dual formulation, i.e. solutions which are real and lead to a positive definite matrix $\vec{H}$. However, some of the real solutions may still be false solutions. These can be pruned even more by substituting the real solutions into the primal formulation and accepting only those with a residual error lower than a given threshold.
    \item \textbf{7r/3s} and \textbf{6r/4s}: For both instances, we first discard the extra receiver (for 7r/3s) or transmitter (for 6r/4s) and call the 6r/3s solver. For each candidate solution of the 6r/3s we trilaterate the extra receiver from $3$ transmitters for 7r/3s, or trilaterate the extra transmitter and offset from $4$ receivers for $6r/4s$. Finally, we substitute the candidate solutions into the primal problem and output the one that leads to the smallest residual error.
    \item \textbf{5r/4s}: As mentioned in the introduction, this problem is interesting because it is not minimal, but it cannot be reduced to minimal configurations either. Using the dual formulation, this configuration has $5+4=9$ unknowns. From the whole compaction matrix, we obtain $2$ independent rank constraints and leaving the last transmitter out, we can obtain $7$ equations in $\vec{H}$,$\vec{b}$ and the first $3$ offsets, having $9$ equations in total. We can thus solve the dual formulation using homotopy continuation and finally trilaterate the last transmitter from $3$ receivers. As this problem is subminimal, the solution is unique. Hence, from the multiple solutions of the dual formulation, we accept the one that leads to the smallest residual error in the primal formulation.
    \item \textbf{5r/5s}: To solve this configuration, we simply leave out the last transmitter and offset, call the 5r/4s solver and finally trilaterate the last transmitter and offset from $4$ receivers. 
\end{itemize}

\section{Results}
\label{ch:results}
\begin{figure*}[t!]
    \centering
    \begin{subfigure}{0.23\linewidth}
        \includegraphics[width=\linewidth]{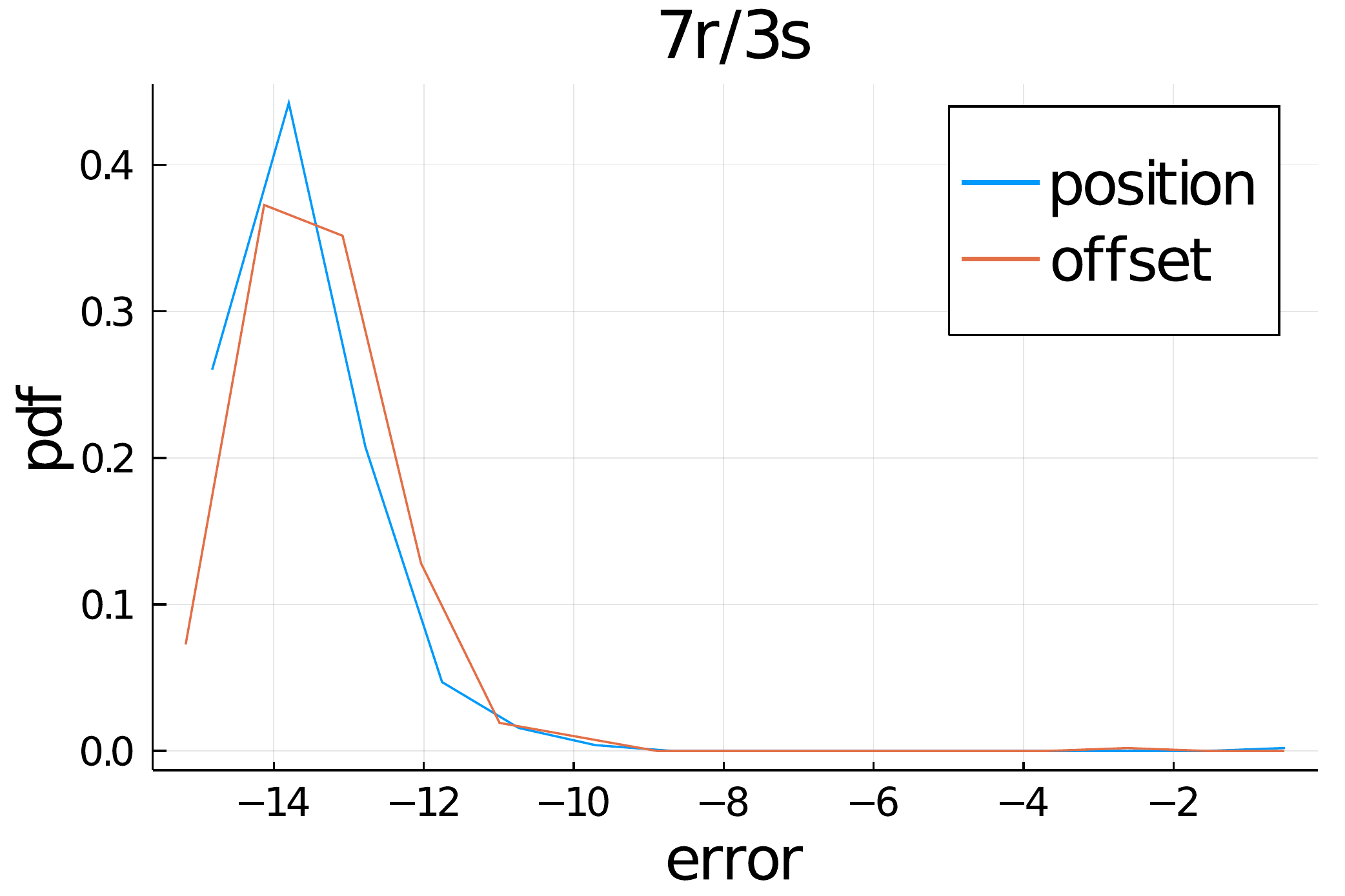}
    \end{subfigure}
    \begin{subfigure}{0.23\textwidth}
        \includegraphics[width=\linewidth]{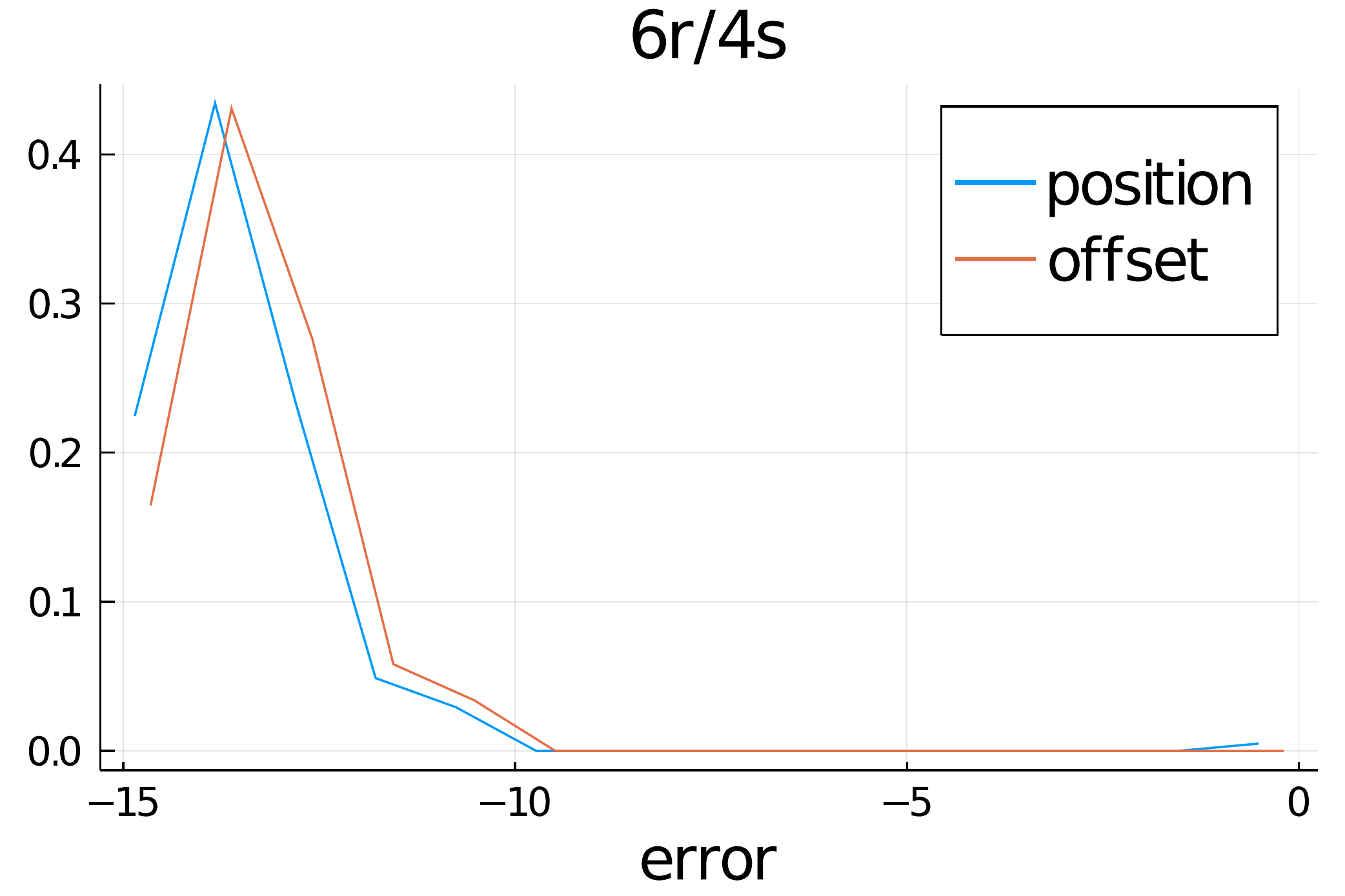}
    \end{subfigure}
    \begin{subfigure}{0.23\textwidth}
        \includegraphics[width=\linewidth]{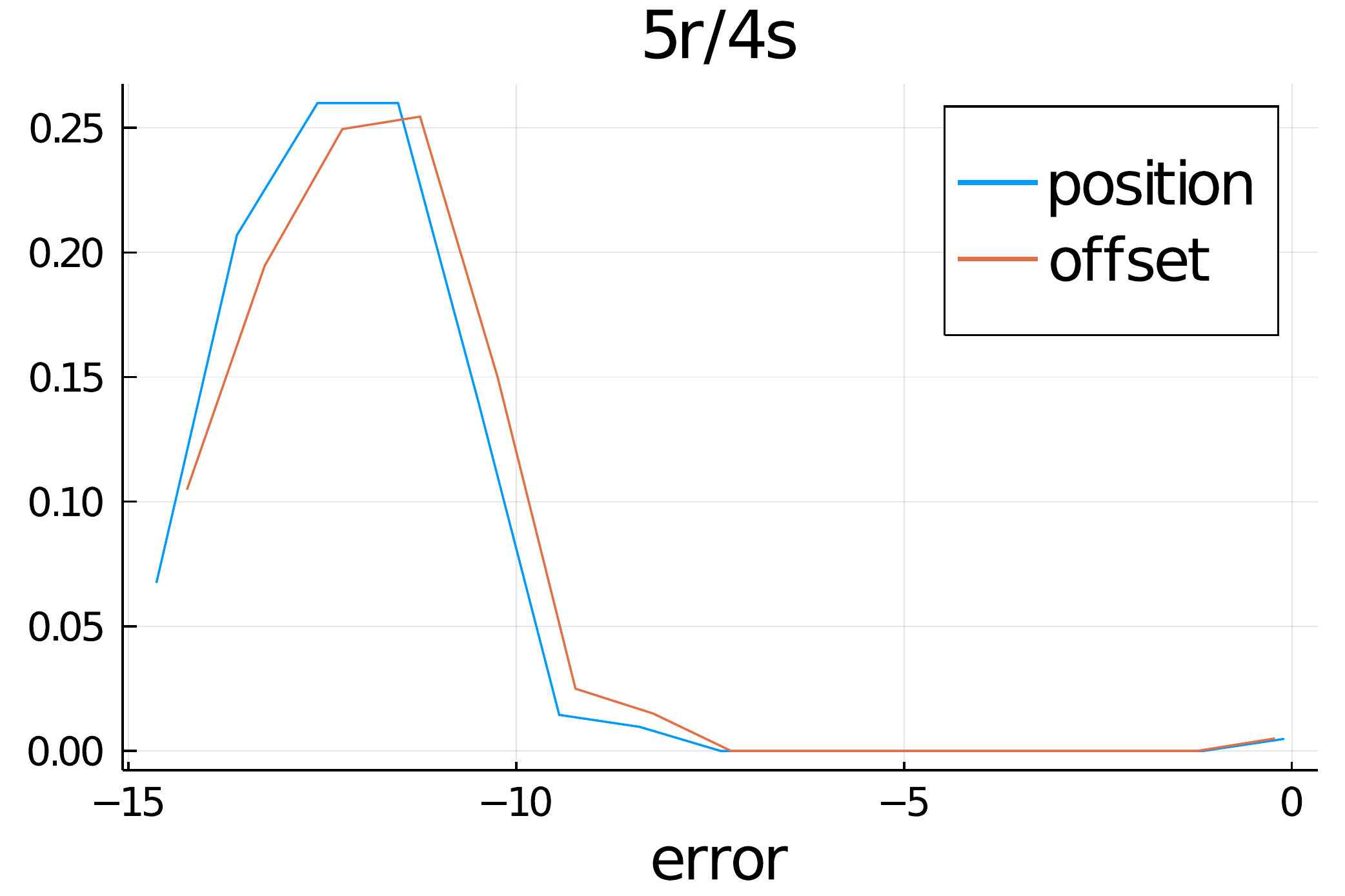}
    \end{subfigure}
    \begin{subfigure}{0.23\textwidth}
        \includegraphics[width=\linewidth]{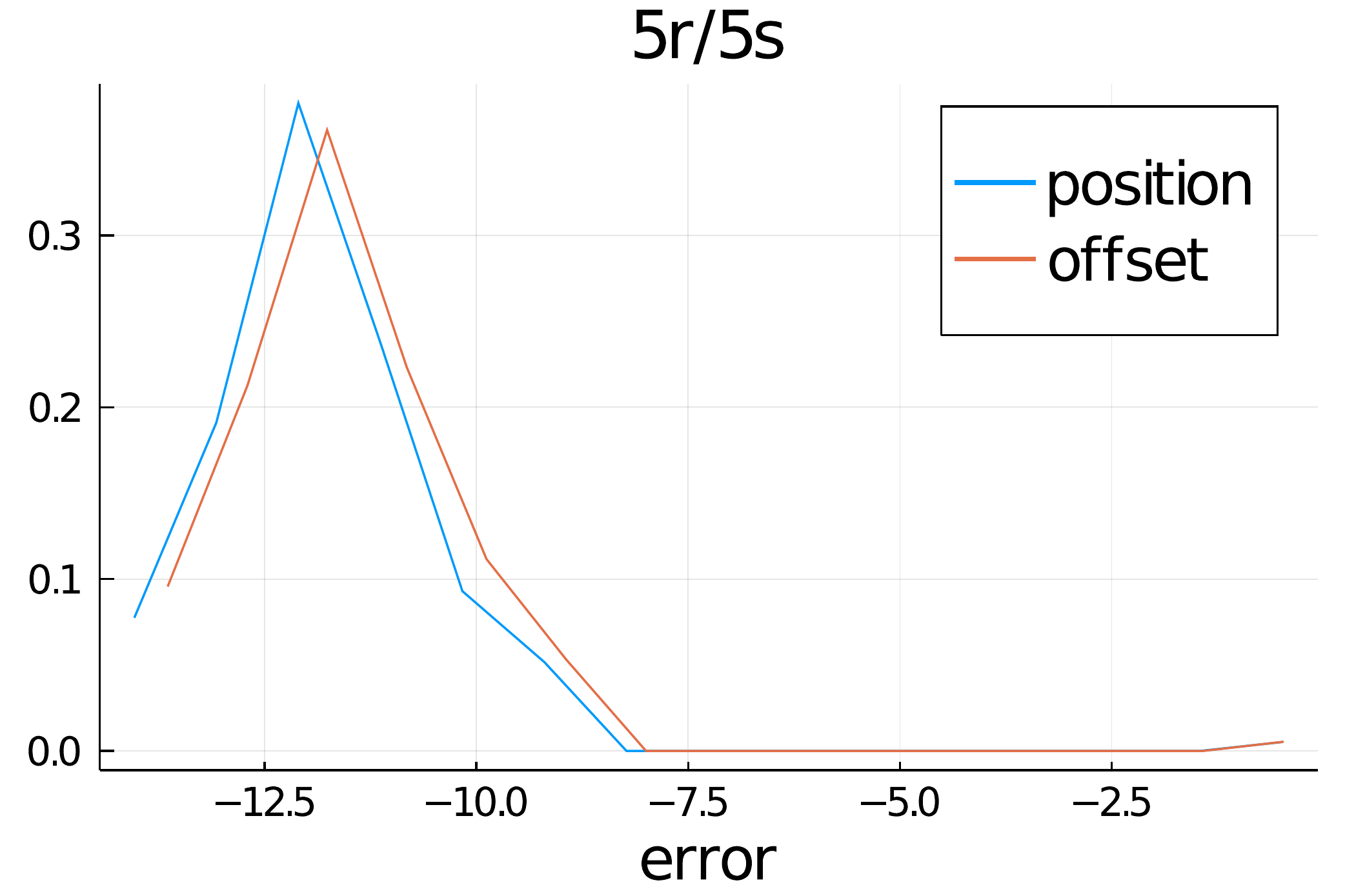}
    \end{subfigure}
    
    \begin{subfigure}{0.23\linewidth}
        \includegraphics[width=\linewidth]{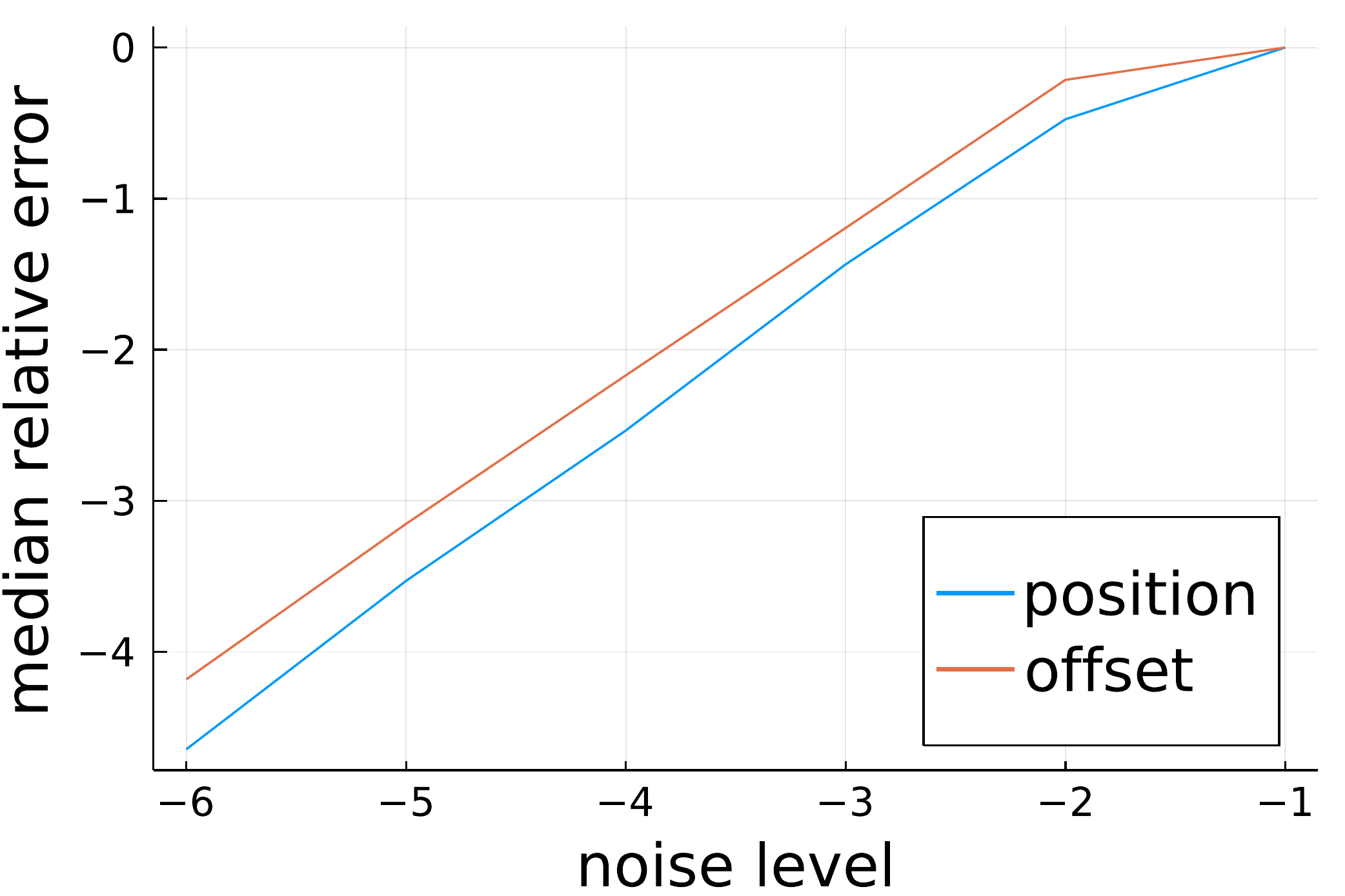}
    \end{subfigure}
    \begin{subfigure}{0.23\textwidth}
        \includegraphics[width=\linewidth]{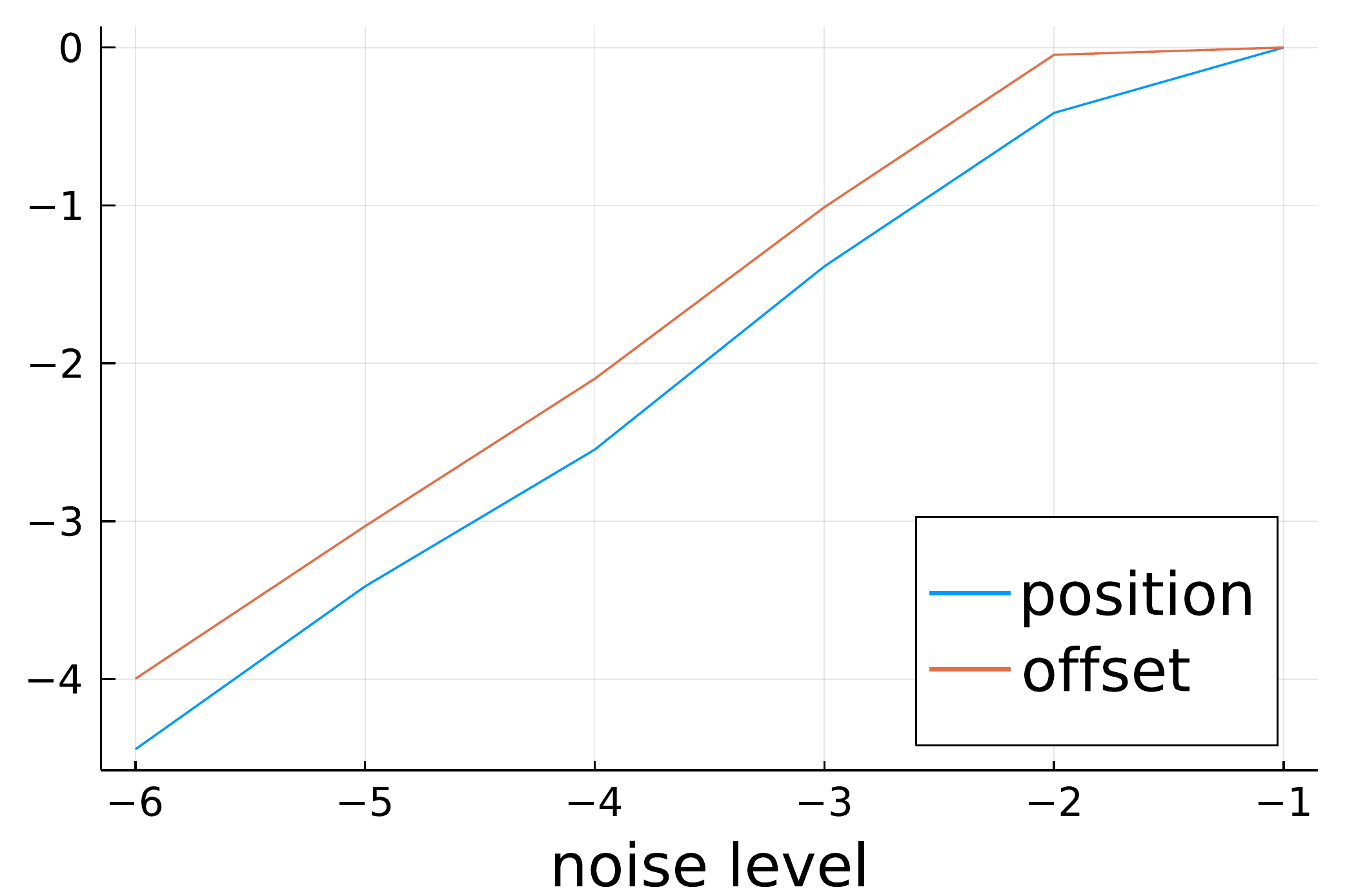}
    \end{subfigure}
    \begin{subfigure}{0.23\textwidth}
        \includegraphics[width=\linewidth]{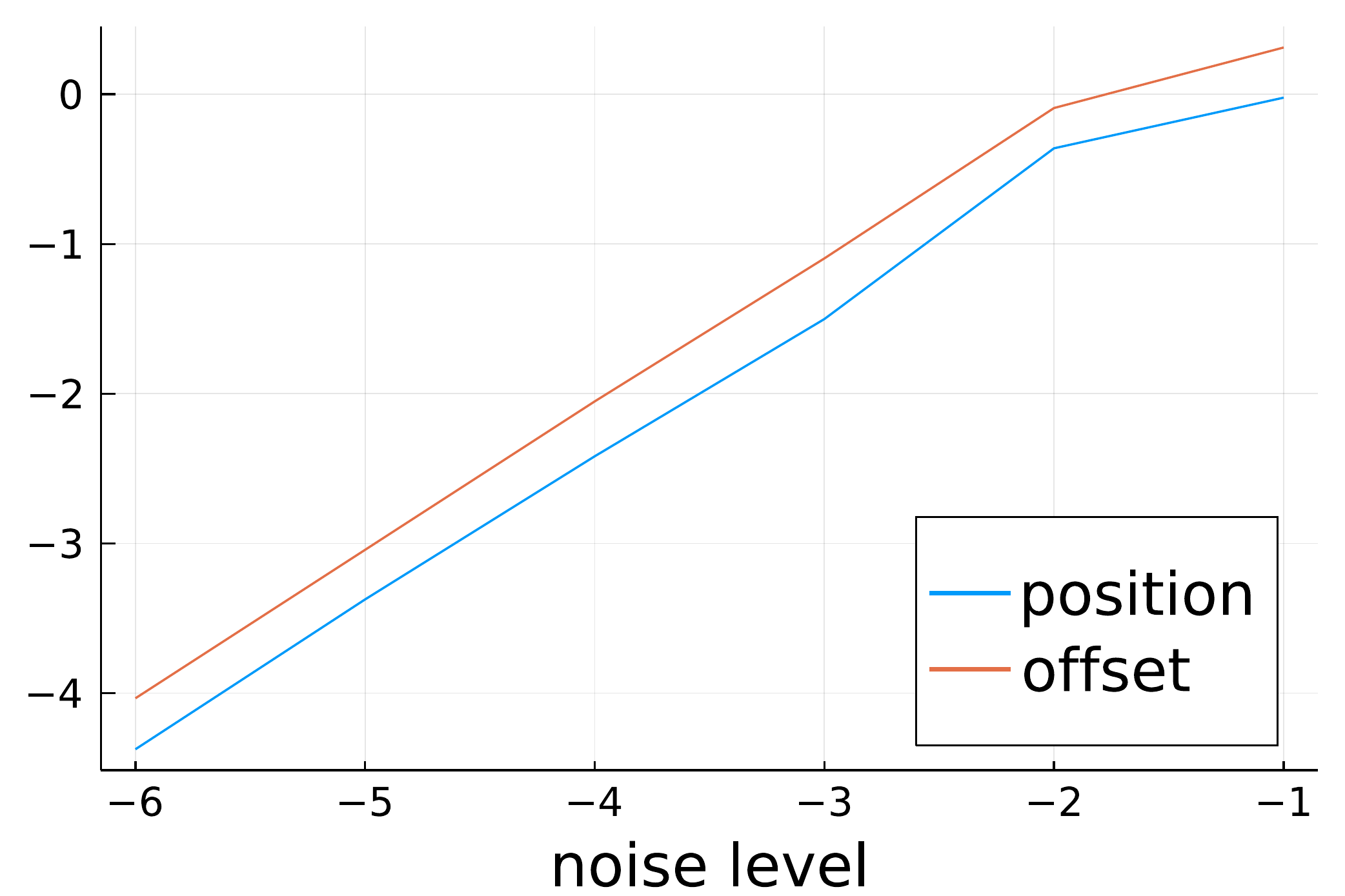}
    \end{subfigure}
    \begin{subfigure}{0.23\textwidth}
        \includegraphics[width=\linewidth]{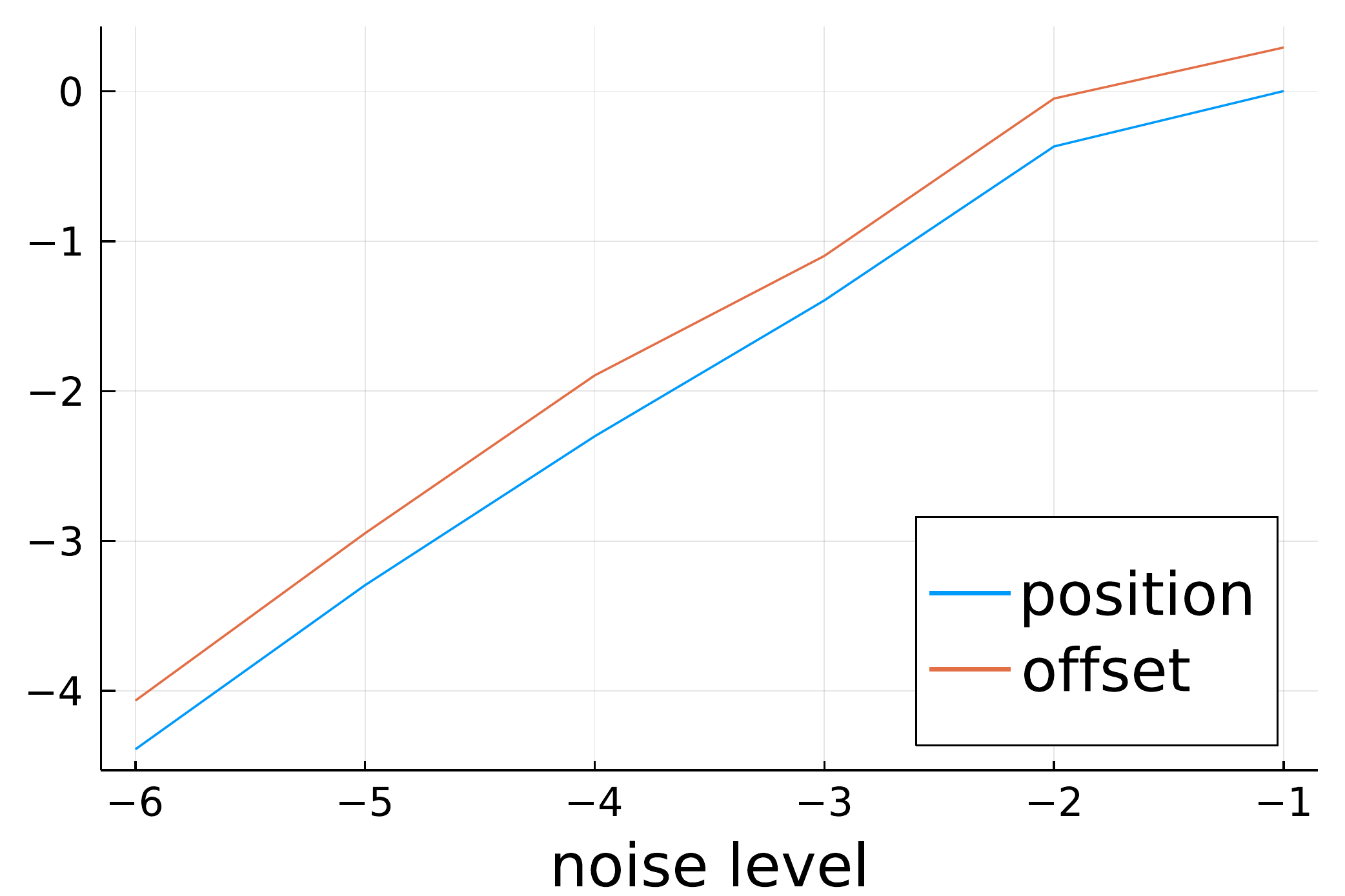}
    \end{subfigure}
    \caption{Quantitative benchmark of our solvers. \textbf{Upper row}: relative error distributions for clean data. \textbf{Lower row}: median relative error at different noise levels, both in logarithmic scale.}
    \label{fig:results}
\end{figure*} 

In this section we discuss benchmarking of the solvers on synthetic data. Th experiments are run with Matlab 2019b and Intel i7-8565U CPU@1.80GHz processor. Each solver ran in approximately 10~s.

\begin{figure}[t]
    \centering
    \begin{subfigure}{0.48\linewidth}
        \includegraphics[width=\linewidth]{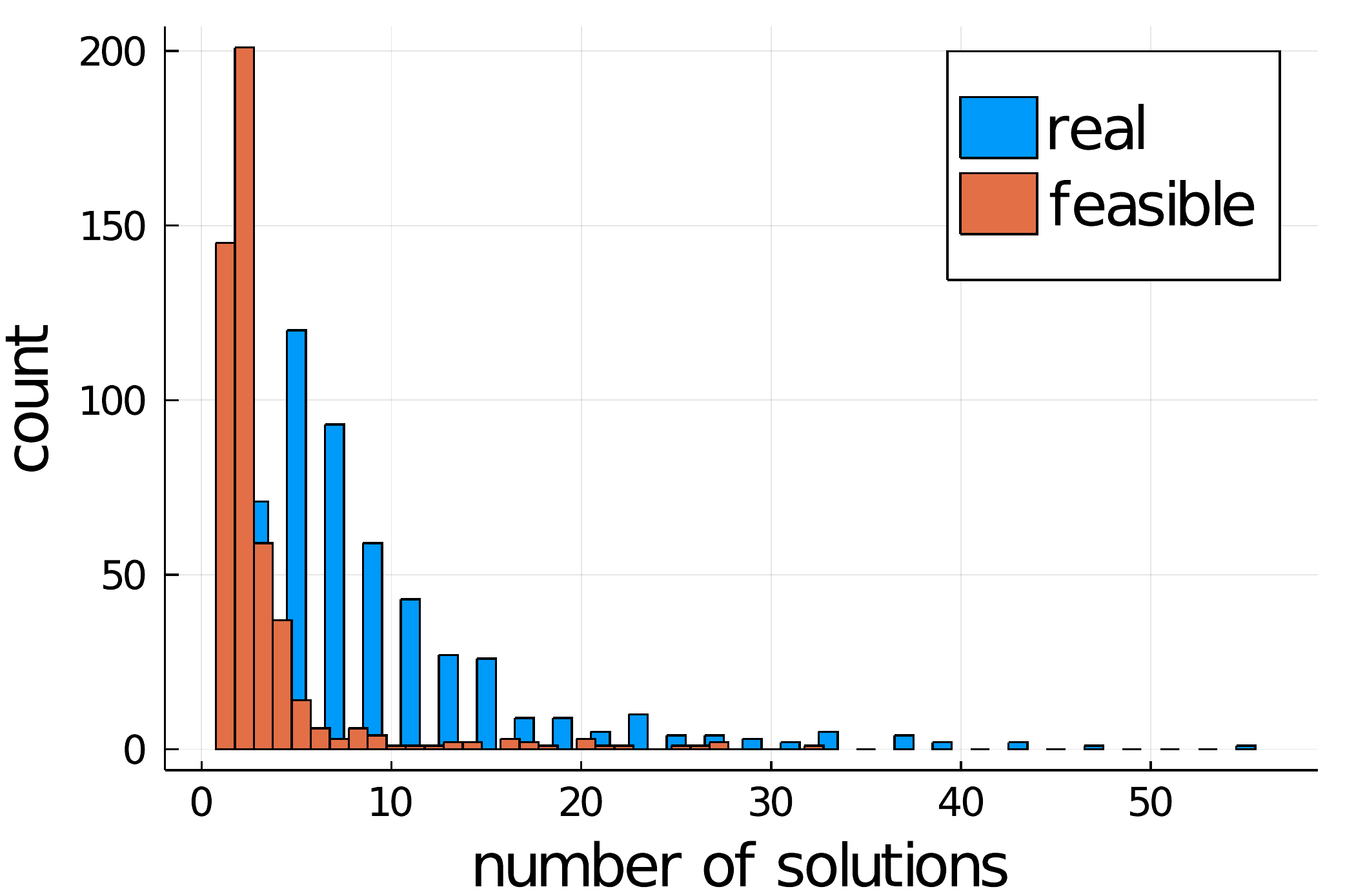}\caption{}\label{fig:hist6r3s}
    \end{subfigure}
    \begin{subfigure}{0.48\linewidth}
        \includegraphics[width=\linewidth]{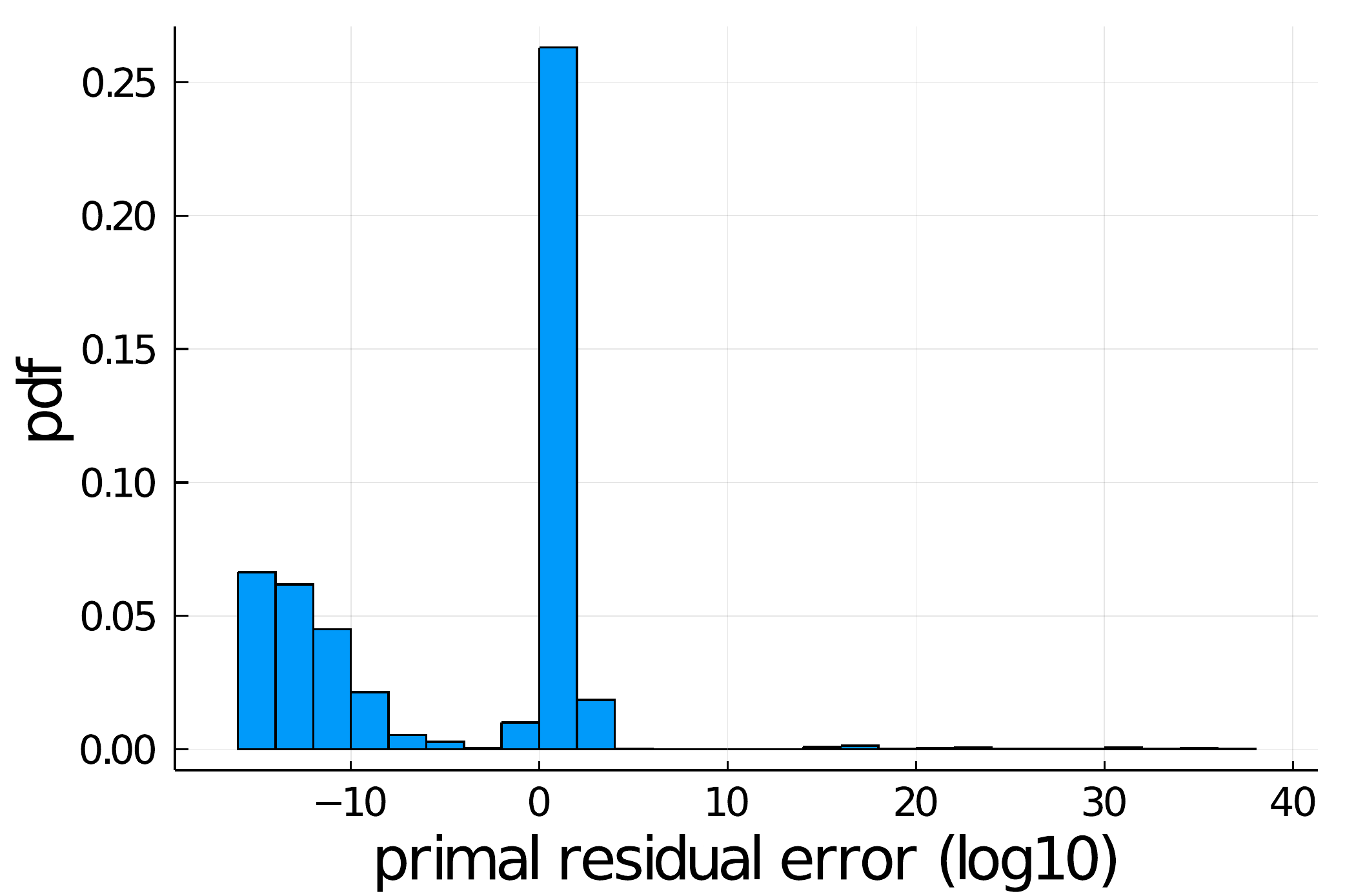}\caption{}\label{fig:error6r3s}
    \end{subfigure}
    \caption{Number of solutions and primal residual errors for 6r/3s solver. }
\end{figure}
\subsection{Minimal solver: 6r3s}
The 6r3s configuration is minimal and as such, it is guaranteed to have a finite, but not unique, number of solutions. It was speculated in \cite{stewenius2005} that the 6r3s configuration would have up to 150 solutions. To investigate this claim, we simulate this configuration with nodes at random positions and offsets and compute the number of real solutions obtained from the dual formulation. The distribution of the number of solutions is depicted in Figure \ref{fig:hist6r3s}. As the histogram and Table~\ref{tab:6r3s} reveal, the dual formulation has from $2$ to $55$ real solutions. However, not all the real solutions of the dual are necessarily solutions of the primal. Indeed, Figure \ref{fig:error6r3s} shows the residual errors when the dual solutions are substituted into the primal system. As this histogram reveals, only some of the dual solutions are solutions of the primal problem and the others are false roots. Practically, this means that we can prune even more the set of real solutions by accepting only those that have a residual error in the primal problem less than a given threshold. In the experiments, we used a threshold of $10^{-10}$ and the distribution of the feasible real solutions is also reported in Figure \ref{fig:hist6r3s} and Table~\ref{tab:6r3s}.
\begin{table}[tbh]
    \centering
    \caption{Number of solutions of the 6r/3s solver}
    \label{tab:6r3s}
    \begin{tabular}{c||c|c|c|c}
         &min&max&mean&st. deviation\\\hline\hline
        real solutions &2&55&8.7&7.6\\
        feas. solutions&1&32&3&3.9
    \end{tabular}
\end{table}
Despite the number of solutions is significantly smaller than what was predicted in \cite{stewenius2005}, still a unique solution cannot be obtained without some extra knowledge of the system. Furthermore, removing the false solutions thresholding the residuals in the dual formulation works fine for clean data. In the presence of noise false roots may also be included, or feasible roots excluded. For this reason, in the next section we also examine subminimal configurations, for which the uniqueness of the solution is guaranteed.

\subsection{Subminimal solvers}
We first benchmark our solvers with clean data, by generating hundreds of random instances of the problems, using a Gaussian distribution for both position and offset. The distributions for position and offset errors are depicted in Figure \ref{fig:results}. As can be noticed, the solvers can achieve very accurate solutions, with the relative error in the order of magnitude $10^{-12}$. We also study how our solvers perform with noisy measurements. As the lower row of Figure \ref{fig:results} shows, our homotopy solver is alone stable also for noisy data even without nonlinear refinement afterwards.

\section{Conclusions}
\label{ch:conclusions}
This paper addressed the problem of sensor networks self-calibration with unsynchronized transmitters. Due to the high number of unknowns and degree of the system, this is a computationally challenging problem to solve. We showed that homotopy continuation offers a powerful tool to overcome these challenges and developed new solvers, stable both for clean and noisy data, which allowed to solve previously unsolved configurations. This opens several perspectives for the applications of homotopy continuation in the domain of network calibration and localization algorithms in general. 

\section*{Acknowledgments}
This work was partially funded by the Academy of Finland project 327912 REPEAT and the Swedish strategic research project ELLIIT.

\bibliographystyle{IEEEtran}
\bibliography{IEEEabrv,IEEEexample}

\end{document}